\documentclass[11pt]{article}
\usepackage{graphicx}
\usepackage{amssymb}
\usepackage{epstopdf}
\input xy
\xyoption{all}

\newtheorem{lemma}{Lemma}[section] 
\newtheorem{propos}[lemma]{Proposition}
\newtheorem{example}[lemma]{Example}
\newtheorem{theorem}[lemma]{Theorem}

\newtheorem{defin}[lemma]{Definition}
\newtheorem{remark}[lemma]{Remark}

\newcommand{\extd}{\mathrm{d}}

\newcommand{\tens}{\mathop{\otimes}}

\newcommand{\id}{\mathrm{id}}
\newcommand{\inc}{\mathrm{inc}}
\newcommand{\im}{\mathrm{im}}
\newcommand{\<}{\langle}
\renewcommand{\>}{\rangle}

\begin{document}
\title{A Leray spectral sequence for noncommutative differential fibrations}
\author{Edwin Beggs\ \dag\ \ \&\ \ Ibtisam Masmali\ \ddag \\ \\
\dag\ College of Science, Swansea University, Wales \\
\ddag\ Jazan University, Saudi Arabia}


\maketitle

\begin{abstract}  This paper describes the Leray spectral sequence associated to a differential fibration. The differential fibration is described by base and total differential graded algebras. The cohomology used is noncommutative differential sheaf cohomology. For this purpose, a sheaf over an algebra is a left module with zero curvature covariant derivative. As a special case, we can recover the Serre spectral sequence for a noncommutative fibration. 
\end{abstract}

\section{Introduction} 

This paper uses the idea of noncommutative sheaf theory introduced in \cite{three}. This is a differential definition, so the algebras involved have to have a differential structure. Essentially
having zero derivative is used to denote `locally constant', which is a term of uncertain meaning for an algebra. Working rather vaguely, one might think of considering the total space of a sheaf over 
a manifold as locally inheriting the differential structure of the manifold, via the homeomorphism between a neighbourhood of a point in the total space and an open set in the base space. This allows us to lift a vector at a point of the base space to a unique vector at every point of the preimage
of that point in the total space. This lifting should allow us to give a covariant derivative on the functions on the total space. Further, the local homeomorphisms suggest that the resulting covariant derivative has zero curvature. In \cite{three} it is shown that a zero curvature covariant derivative on a module really does allow us to reproduce some of the main results of sheaf cohomology. 
In this paper we shall consider another of the main results of sheaf cohomology, the Leray
spectral sequence.

Ideally it would be nice to have a definition which did not involve differential structures, but there are several comments to be made on this: 
When Connes calculated the cyclic cohomology of the noncommutative torus
\cite{ConnesIHES}, he used a subalgebra of rapidly decreasing sequences, effectively 
placing differential methods at the heart of noncommutative cohomology. It is not obvious what a \textit{calculable} purely algebraic (probably read $C^*$ algebraic) sheaf cohomology theory would be -- though maybe the theory of quantales \cite{MulQuant} might give a clue.  
Secondly, even if there were a non-differential definition, it would likely be complementary to the differential definition. The relation between de Rham and topological cohomology theories is fundamental to a lot of mathematics, it would make no sense to delete either. Finally, in
mathematics today, differential graded algebras arising from 
several constructions are considered interesting objects in their own right, and many applications to Physics are phrased in terms of differential forms or vector fields. 

There are four main motivations behind this paper. One is that the Leray spectral sequence seems a natural continuation from the sheaf theory and Serre spectral sequence in \cite{three}. Another is a step in finding an analogue of the Borel-Weil-Bott theorem for representations of quantum groups
(see \cite{BWquant93}). One motivation we should look at in more detail is contained in the papers \cite{25,26}. These papers are about noncommutative fibrations. The differences in approach can be summarised in two sentences: We require that the algebras have differential structures, and \cite{25,26} do not. The papers \cite{25,26} require that the base is commutative, and we do not. One interesting point is that the method of \cite{26} makes use of the classical Leray spectral sequence
 of a fibration with base a simplicial complex. The fourth motivation is noncommutative algebraic topology, where we would define a fibration on a category whose objects were differential graded algebras. The interesting question is then whether there is a corresponding idea of cofibration in the sense of model categories \cite{quillModel}. 
 
 The example of the noncommutative Hopf fibration in \cite{three} shows that a differential fibration need not have a commutative base.
 The example in Section \ref{se1} was made by taking a differential picture of a fibration given as an example in \cite{25} (the base is the functions on the circle), and so it can be considered a noncommutative fibration in both senses. It would be useful to consider whether higher dimensional constructions, such as the  4-dimensional orthogonal quantum sphere in \cite{48}, also give examples of differential fibrations. As differential calculi on finite groups are quite well understood (e.g.\ see \cite{17,Ma:rief}), it would be interesting to ask what a differential fibration corresponds to in this context. From the point of view of methods in mathematical Physics, the quantisation of twistor theory
 (see \cite{BraMa}) is likely to provide some examples. 

This paper is based on part of the content of the Ph.D.\ thesis \cite{MasThesis}.

\section{Spectral sequences}
This is standard material, and we use  \cite{11} as a reference. We will give quite general definitions, but likely not the most general possible. 

\subsection{What is a spectral sequence?}
A spectral sequence consists of series of pages (indexed by $r$) and objects $ \mathcal{E}^{p,q}_{r}$ (e.g.\ vector spaces), where $r,p,q$ are integers. We take $r\geq 1$ and $p,q \geq 0$ , and set $ \mathcal{E}^{p,q}_{r} = 0$ if $p < 0$ or  $q < 0$ . There is a differential $$\extd_{r} :  \mathcal{E}^{p,q}_{r} \longrightarrow  \mathcal{E}^{p+r,q+1-r}_{r}$$
such that $\extd_{r}\extd_{r} = 0$.
As $\extd_{r}\extd_{r} = 0$, we can take a quotient (in our case, quotient of vector spaces) $$\frac{ \ker  \, \extd_{r} :  \mathcal{E}^{p,q}_{r} \rightarrow  \mathcal{E}^{p+r,q+1-r}_{r}}{\im \, \extd_{r} :  \mathcal{E}^{p-r,q+r-1}_{r} \rightarrow  \mathcal{E}^{p,q}_{r}} = H^{p,q}_{r}$$ Then the rule for going from page $r$ to page $r+1$ is $ \mathcal{E}^{p,q}_{r+1} = H^{p,q}_{r}$.
The maps $d_{r+1}$ are given by a detailed formula on $H^{p,q}_{r}$.
The idea is that eventually the $ \mathcal{E}^{p,q}_{r}$ will become fixed for $r$ large enough.
The spectral sequence is said to converge to these limiting cases $ \mathcal{E}^{p,q}_{\infty}$ as $r$ increases. 

\subsection{The spectral sequence of filtration}\label{a53}

A decreasing filtration of a vector space $V$ is a sequence of subspaces $F^m V$
($m\in\mathbb{N}$) for which $F^{m+1}V \subset F^m V$.
The reader should refer to \cite{11} for the details of the homological algebra used to construct the spectral sequence. We will merely quote the results. 

\begin{remark}\label{sprem} Start with a differential graded module $C^n$ (for $n\ge 0$) and $ \extd :C^n \to C^{n+1}$ with $ \extd^2=0$. Suppose that $C$ has a filtration  $F^m C\subset C=\oplus_{n\ge 0}C^n$ for $m\ge 0$ so that:\\
 (1)\quad $ \extd F^m C \subset F^m C$ for all $m\ge 0$ (i.e.\ the filtration is preserved by $ \extd$); \\
 (2)\quad $F^{m+1} C\subset F^m C$ for all $m\ge 0$ (i.e.\ the filtration is decreasing); \\
 (3)\quad $F^0 C=C$ and $F^m C^n=F^m C\cap C^n=\{0\}$ for all $m>n$ (a boundedness condition). \\
 Then there is a spectral sequence $(\mathcal{E}_r^{*,*},  \extd_r)$ for $r\ge 1$ ($r$ counts the page of the spectral sequence) with $ \extd_r$ of bidegree $(r,1-r)$ and
 \begin{eqnarray}\label{b7}
 \mathcal{E}_1^{p,q} &=& H^{p+q}(F^pC/F^{p+1}C) \cr
&=& \frac{{\rm ker}\, \extd:F^pC^{p+q}/F^{p+1}C^{p+q}\to F^pC^{p+q+1}/F^{p+1}C^{p+q+1}}{{\rm im}\, \extd:F^pC^{p+q-1}/F^{p+1}C^{p+q-1}\to F^pC^{p+q}/F^{p+1}C^{p+q}}\ .
 \end{eqnarray}
 In more detail, we define
 \begin{eqnarray*}
 Z_{r}^{p,q} &=& F^{p} C^{p+q} \cap  \extd^{-1}(F^{p+r} C^{p+q+1})\ ,\cr
 B_{r}^{p,q} &=& F^{p} C^{p+q} \cap \extd(F^{p-r} C^{p+q-1})\ ,\cr
 \mathcal{E}_{r}^{p,q} &=& Z_r^{p,q}/(Z_{r-1}^{p+1,q-1}+B_{r-1}^{p,q})\ 
 .\end{eqnarray*}
 The differential $ \extd_{r}:\mathcal{E}_{r}^{p,q} \to \mathcal{E}_{r}^{p+r,q-r+1}$ is the map induced on quotienting $ \extd:Z_{r}^{p,q} \to Z_{r}^{p+r,q-r+1}$. The diligent reader should remember an important point here, when reading the seemingly innumerable differentials in the pages to come. There is really only one differential $\extd$ -- its domain or codomain may be different subspaces with different quotients applied, but the same $\extd$ lies behind them all. 
 
  The spectral sequence converges to $H^*(C, \extd)$ in the sense that
 \begin{eqnarray*}\mathcal{E}_\infty^{p,q} \cong \frac{F^p H^{p+q}(C, \extd)}{F^{p+1}H^{p+q}(C, \extd)}\ ,
 \end{eqnarray*}
 where $F^p H^*(C, \extd)$ is the image of the map $H^*(F^p C, \extd)\to H^*(C, \extd)$ induced by inclusion $F^p C\to C$.
 \end{remark}

\subsection{The classical Leray spectral sequence}
The statement of the general Leray spectral sequence can be found in \cite{28}. We shall omit the supports and the subsets as we are only currently interested in a non commutative analogue of the spectral sequence.
Then the statement reads that, given $f : X \rightarrow Y$ and $\mathcal{S}$ a sheaf on $X$, that there is a spectral sequence$$E^{pq}_{2} = H^{p}(Y, H^{q}(f,f \vert \mathcal{S}))$$ converging to $H^{p+q}(X,\mathcal{S})$.
Here $H^{q}(f,f \vert \mathcal{S}) $ is a sheaf on $Y$ which is given by the presheaf for an open $U \subset Y$  $$U\longmapsto H^{q}(f^{-1}U; \mathcal{S} \vert_{f^{-1}U}).$$
Here $f^{-1}U$ is an open set of $X$, and $\mathcal{S} \vert_{f^{-1}U}$ is the sheaf $\mathcal{S}$ restricted to this open set.

We shall consider the special case of a differential fibration. This is the background to the Serre spectral sequence, but we consider a sheaf on the total space.
The Leray spectral sequence of a fibration is a spectral sequence whose input is the cohomology of the base space $B$ with coefficients in the cohomology of the fiber $F$, and converges to the cohomology of the total space $E$. Here$$\pi : E \rightarrow B$$ is a fibration with fiber $F$. The difference of this from the Serre spectral sequence is that the cohomology may have coefficients in a sheaf on $E$.

\section{Noncommutative differential calculi and sheaf theory}
Take a possibly noncommutative algebra  $A$. Then a differential calculus $(\Omega^*A,\extd)$
is given by the following. 

\begin{defin}\label{anwar}
A differential calculus $(\Omega^*A,\extd)$ on $A$
 consists of vector spaces $\Omega^{n}A$ with operators  $\wedge$ and $\extd$ so that  \\
1) $\wedge : \Omega^r A \otimes \Omega^m A \longrightarrow \Omega^{r+m} A$  is associative (we do not assume any graded commutative property) \\
2) $\Omega^0 A = A $ \\
3) $\extd : \Omega^n A \rightarrow \Omega^{n+1}A $ with $\extd^2 =0$ \\
4) $\extd(\xi \wedge \eta ) = \extd\xi \wedge \eta + (-1)^r \xi \wedge \extd\eta$  for $\xi \in \Omega^r A$  \\
5) $\Omega^1 A \wedge \Omega^n A = \Omega^{n+1} A$ . \\
6) $A.\extd A = \Omega^{1}A$ 
\end{defin}

Note that many differential graded algebras do not obey (5), but those in classical differential geometry do, and it will be true in all our examples. There is only one place where we require (5), and we will point it out at the time.

A special case of $\wedge$
shows that each $\Omega^n A$ is an $A$-bimodule. 
We will often use $\vert \xi \vert$ for the degree of $ \xi $, if $\xi \in \Omega^{n}A$, then $\vert \xi \vert = n$.

In the differential graded $(\Omega^{n} A,\wedge,\extd)$, we have $\extd^{2} = 0$. This means that 
$$\im \,  \extd : \Omega^{n-1} A \longrightarrow \Omega^{n} A  \subset  \ker \, \extd : \Omega^{n} A \longrightarrow \Omega^{n+1} A\ .$$
Then we define the de Rham cohomology as
 $$H^{n}_{\mathrm{dR}}(A)\ =\ \frac{\ker  \,  \extd : \Omega^{n} A \longrightarrow \Omega^{n+1} A}{\im \,  \extd : \Omega^{n-1} A \longrightarrow \Omega^{n} A}\ .$$
We give the usual idea of covariant derivatives on left $A$ modules by using the left Liebnitz rule:

\begin{defin}\label{de1}
Given a left  A-module E , a left  A-covariant  derivative is a map $\nabla : E \rightarrow \Omega^1 A \otimes_{A} E$ which obeys the condition $\nabla ( a.e) = \mathrm{d}a \otimes e + a . \nabla e $  for all $e \in E$ and $a \in A$.
\end{defin}

After the fashion of the de-Rham complex, we can attempt to extend the covariant derivative to a complex as follows:

\begin{defin}\cite{three}
Given $(E,\nabla )$ a left $A$-module with covariant derivative, define $$ \nabla^{[n]}  : \Omega^{n} A \otimes_{A} E \rightarrow \Omega^{n+1} A \otimes_{A} E ,   \quad   \omega \otimes e \mapsto \extd\omega \otimes e + (-1)^{n} \omega \wedge \nabla e .$$ Then the curvature is defined as $R = \nabla^{[1]} \nabla : E \rightarrow \Omega^{2} A \otimes E$, and is a left A-module map. 
The covariant derivative is called flat if the curvature is zero.
\end{defin} 

However, the curvature forms an obstruction to setting up a cohomology, as we now show:

\begin{propos}\cite{three}
For all $n \geq 0$, $\nabla^{[n+1]} \circ \nabla^{[n]} = \id \wedge R : \Omega^{n} A \otimes_{A} E \rightarrow \Omega^{n+2} A \otimes_{A} E .$
\end{propos}

 We can now use this in a definition of a noncommutative sheaf \cite{three}.

\begin{defin}\label{key49}\cite{three}
Given $(E,\nabla )$ a left $A$-module with covariant derivative and zero curvature, define $H^*(A ; E, \nabla )$ to be the cohomology of the cochain complex $$E \stackrel{\nabla} \longrightarrow \Omega^{1} A \otimes_{A}E
\stackrel{\nabla^{[1]}} \longrightarrow \Omega^{2} A \otimes_{A} E \stackrel{\nabla^{[2]}}\longrightarrow ........$$ Note that $ H^{0}(E , \nabla ) =  \{ e \in E : \nabla e = 0 \}$, the flat sections of $E$. We will often write $H^* (A ;E)$ where there is no danger of confusing the covariant derivative .
\end{defin}

We will take this opportunity to make a couple of well known 
 statements about modules over algebras which we will use, as it may make the reading later easier for non-experts (see e.g.\ \cite{44}) .

\begin{defin} A right A-module $E$ is flat if every short exact sequence of left A-modules $$0 \longrightarrow L \longrightarrow M \longrightarrow N  \longrightarrow 0$$
gives another short  exact sequence $$0 \longrightarrow E\otimes_{A} L \longrightarrow E\otimes_{A} M \longrightarrow E\otimes_{A} N \longrightarrow  0.$$
Similarly, a left A-module $F$ is called flat if $-\tens_A F$ preserves exactness
of short sequences of right modules. 
\end{defin}

\begin{lemma}\label{raneem}
Given two short exact sequences of modules (left or right),
\begin{eqnarray*}
 && 0 \longrightarrow U \stackrel{t}\longrightarrow V \stackrel{f}\longrightarrow W \longrightarrow 0 \ ,\cr
&& 0 \longrightarrow U \stackrel{t}\longrightarrow V \stackrel{g}\longrightarrow X \longrightarrow 0\ ,
\end{eqnarray*}
there is an isomorphism $h : W \longrightarrow X $ given by $h(w) = g(v)$, where $f(v) = w$. 
\end{lemma}

\section{Differential fibrations and the Serre spectral sequence}

\subsection{A simple differential fibration}\label{a52}
The reader may take this section as a justification of why the definition of a noncommutative
differential fibration which we will give in Definition \ref{b61} is reasonable. 
Take a trivial fibration $\pi:\mathbb{R}^n \times \mathbb{R}^m \to \mathbb{R}^n$ given by
$$(x_{1} , . . . .,x_{n},y_{1} , . . . .,y_{m}) \longmapsto (x_{1} , . . . .,x_{n})\ .$$
Here the base space is $B = \mathbb{R}^n$, the fiber is $\mathbb{R}^m$, and the total space is $E = \mathbb{R}^{n+m}$.
We can write a basis for the differential forms on the total space, putting the $B$ terms (the $\extd x_{i}$) first. A form of degree $p$ in the base and $q$ in the fiber (total degree $p+q$) is
$$\extd x_{\iota_{1}} \wedge . . . .\wedge \extd x_{\iota_{p}} \wedge \extd y_{j_{1}} \wedge . . . .\wedge \extd y_{j_{q}}\ ,$$
for example $\extd x_{2} \wedge \extd x_{4} \wedge \extd y_{1} \wedge \extd y_{7} \wedge \extd y_{9}$, 
If we have the projection map $\pi : E \longrightarrow B$, we can write our example form as
$$\alpha \ =\ \pi^*(\extd x_{2} \wedge \extd x_{4}) \wedge (\extd y_{1} \wedge \extd y_{7} \wedge \extd y_{9})$$
so we have a form in $\pi^* \Omega^2 B \wedge \Omega^3 E$. Another element of $\pi^* \Omega^2 B \wedge \Omega^3 E$ might be
$$\beta\ =\ \pi^*(\extd x_{2} \wedge \extd x_{4}) \wedge (\extd x_{3} \wedge \extd y_{1} \wedge \extd y_{_{7}}).$$
Note, we now just look at $\Omega^3 E$, not the forms in the fiber direction, as in the noncommutative case we will not know (at least in the begining) what the fiber is. We need to describe  the forms on the fiber space more indirectly.
Now look at the vector space quotient
\begin{eqnarray}\label{vcfhgmk}
\frac{\pi^* \Omega^2 B \wedge \Omega^3 E}{\pi^* \Omega^3 B \wedge \Omega^2 E}\ .
\end{eqnarray}
Here $\beta$ is also an element of the bottom line of (\ref{vcfhgmk}), as we could write
$$\beta = \pi^*(\extd x_{2} \wedge \extd x_{4} \wedge \extd x_{3}) \wedge (\extd y_{1} \wedge \extd y_{_{7}})$$
so, denoting the quotient by square brackets, $[\beta] = 0$. On the other hand, $\alpha$ is not in the bottom line of (\ref{vcfhgmk}), so $[\alpha] \neq 0$. We can now use $$\frac{\pi^* \Omega^p B \wedge \Omega^q E}{\pi^* \Omega^{p+1} B \wedge \Omega^{q-1} E}$$
to denote the forms on the total space which are of degree $p$ in the base and degree $q$ in the fiber, without explicitly having any coordinates for the fiber. This is just the idea of a 
noncommutative differential fibration. 

\subsection{Noncommutative differential fibrations}
In Subsection \ref{a52} we had a topological fibration $\pi:\mathbb{R}^{m+n}\to \mathbb{R}^n$.
For algebras, we will reverse the arrows, and look at $ \iota : B \rightarrow A$, where $B$ is the `base algebra' and $A$ is the `total algebra'. 

Suppose that both $A$ and $B$ have differential calculi, and that the algebra map $\iota : B \rightarrow A$
is differentiable. This means that $\iota : B \rightarrow A$ extends to a map of differential graded algebras
$\iota_* : \Omega^*B \rightarrow \Omega^*A$, and in particular that $\extd\,\iota_*=\iota_*\,\extd$
and $\iota_*\,\wedge=\wedge\,(\iota_*\tens \iota_*)$. 
Now we set
\begin{eqnarray} \label{cvhgsuv}
D_{p,q} = \iota_{*} \Omega^{p} B \wedge \Omega^{q} A\quad\mathrm{and} \quad
N_{p,q} = \frac{D_{p,q}}{D_{p+1,q-1}}\ ,\quad N_{p,0} \,=\, \iota_*\Omega^{p} B.A\ .
\end{eqnarray}
Now we can finally define a differential fibration, remembering that we use $[\ ]$ to denote equivalence class in the quotient in (\ref{cvhgsuv}):

\begin{defin}\label{b61}
$\iota : B \longrightarrow A$ is a differential fibration if the map$$\xi \otimes [x] \longrightarrow [\iota_{*} \xi \wedge x]$$ gives an isomorphism from $\Omega^{p}B \otimes_{B} N_{0,q}$ to $N_{p,q}$ for all $p,q\ge 0$.
\end{defin}

\begin{example}\label{f2}
(See section 8.5 of \cite{three}.) Given the left covariant calculus on the quantum group $SU_{q}(2)$ given by Woronowicz \cite{worondiff}, the corresponding differential calculus on the quantum sphere $S^{2}_{q}$ gives a differential fibration $$\iota : S^{2}_{q} \longrightarrow SU_{q}(2)\ .$$Here the algebra $SU_{q}^{2}$ is the invariants of $SU_{q}(2)$ under a circle action, and $\iota$ is just the inclusion.
\end{example}

\noindent We will give another example in Section \ref{se1}. Now we have the following version of the Serre spectral sequence from \cite{three}. 

\begin{theorem}
Suppose that $\iota : B \rightarrow X$ is a differential fibration. Then there is a spectral sequence converging to $H^{*}_{dR}(A)$ with $$E^{p,q}_{2} \cong H^{p}(B ; H^{q}(N_{0,*} ), \nabla )\ .$$
\end{theorem}

Here $\nabla$ is a zero curvature covariant derivative on the left $B$-modules $N_{0,n}$, whose construction we will not go further into, as we are about to something more general.

\section{The noncommutative Leray spectral sequence}

 \subsection{A filtration of a cochain complex}\label{a54}
 We suppose that $E$ is a left A module, with a left covariant derivative $$\nabla : E \longrightarrow \Omega^1 A \otimes_{A}E$$
 and that this covariant derivative is flat, i.e.\ that its curvature vanishes. Then $\nabla^{[n]} : \Omega^{n}A\otimes_{A}E \longrightarrow \Omega^{n+1}A \otimes_{A}E$ is a cochain complex (see definition \ref{key49}). 
 Suppose that $\iota_* : \Omega^*B \longrightarrow \Omega^*A$ is a map of differential graded algebras. We define a filtration of $\Omega^{n}A \otimes_{A}E$ by
\begin{eqnarray}\label{rana}
F^m(\Omega^{n}A \otimes_{A} E) = \left\{ \begin{array}{ll}
         \iota_{*} \Omega^m B \wedge \Omega^{n-m} A \otimes_{A} E & \mbox{ $0 \leq m \leq n$};\\
        0 & \mathrm{otherwise}.\end{array} \right. 
\end{eqnarray} 
\begin{propos}
The filtration in (\ref{rana}) satisfies the conditions of remark \ref{sprem}.\\
\textbf{Proof:}
 First  \quad $ F^0 (\Omega^n A \otimes_{A} E) = \iota_{*}\Omega^0 B \wedge \Omega^n A \otimes_{A} E,$
 but $1 \in \iota_{*} \Omega^0 B = \iota_{*}B$,   so  $ F^0 (\Omega^n A \otimes_{A} E) =  \Omega^n A \otimes_{A} E.$.
 
 To show it is decreasing, (using condition (5) from definition \ref{anwar})
 \begin{eqnarray*}
  F^{m+1} (\Omega^n A \otimes_{A} E) &=& \iota_{*}\Omega^{m+1} B \wedge \Omega^{n-m-1} A \otimes_{A} E \cr
  &=& \iota_{*} \Omega^m B \wedge (\iota_{*} \Omega^{1} B \wedge \Omega^{n-m-1}A) \otimes_{A} E\cr
  &\subset &  \iota_{*} \Omega^m B \wedge \Omega^{n-m}A \otimes_{A} E \cr
  &\subset & F^m (\Omega^n A \otimes_{A} E)\ .
\end{eqnarray*}
To show that the filtration is preserved by $\extd$, take $\iota_{*}
 \xi \wedge \eta \otimes e \in F^m (\Omega^n A \otimes_{A} E)$
where $\xi \in \Omega^m B$, and  $\eta \in \Omega^{n-m} A$. Then 
$$\extd(\iota_{*}  \xi \wedge \eta \otimes e) = \iota_{*}  \extd\xi \wedge \eta \otimes e + (-1)^m \iota_{*}  \xi \wedge \extd\eta \otimes e +(-1)^n \iota_{*}  \xi \wedge \eta \wedge \nabla e$$
This is in $F^m C$, as the first term is in $F^{m+1}C \subset F^m C$, and the other two are in $F^m C$.\quad $\square$
\end{propos}
Now we have a spectral sequence which converges to $H^{*}_{d R}(A ; E)$. All we have to do is to find the first and second pages of the spectral sequence, though this is quite lengthy.

\subsection{The first page of the spectral sequence}\label{first}
From section \ref{a53}, to use the filtration in section \ref{a54} we need to work with
\begin{eqnarray}  \label{cbdhsiouv}
M_{p,q} = \frac{F^p C^{p+q}}{F^{p+1} C^{p+q}} = \frac{\iota_{*}\Omega^p B \wedge \Omega^q A \otimes_{A} E}{\iota_{*}\Omega^{p+1} B \wedge \Omega^{q-1} A \otimes_{A} E}
\end{eqnarray}
Then we look, for $p$ fixed (following (\ref{b7})), at the sequence 
\begin{eqnarray}\label{b60}
\cdots M_{p,q-1} \stackrel{\extd} \longrightarrow M_{p,q} \stackrel{\extd} \longrightarrow M_{p,q+1} \stackrel{\extd} \longrightarrow  \cdots
\end{eqnarray}
as the cohomology of this sequence gives the first page of the spectral sequence.
Denote the quotient in $M_{p,q}$ by $[ \quad ]_{p,q}$, so if $x \in \iota_{*} \Omega^{p}B \wedge \Omega^{q} A \otimes_{A} E$, then $[x]_{p,q} \in M_{p,q}$.
Then we have a map of left $B$ modules $$\Omega^{p}B \otimes_{B} M_{0,q} \longrightarrow M_{p,q}\ ,\quad \xi \otimes [y]_{0.q} \longmapsto [\iota_{* }\xi \wedge y]_{pq}.$$
Here $y \in \Omega^{q} A \otimes_{A} E$ and the left action of $b \in B$ on $y$ is $\iota(b)y$.

\begin{propos}\label{aa66}
If $E$ is flat as a left A module, then $N_{p,q} \otimes_{A} E \cong M_{p,q}$ with isomorphism $[z] \otimes e \longmapsto [z \otimes e]_{p,q}$.\end{propos}
\textbf{Proof:} We have, by definition, a short exact sequence
using notation from (\ref{cvhgsuv}), where $\inc$ is inclusion and $[\, \,]$ is quotient, $$0 \longrightarrow D_{p+1,q-1} \stackrel{\mathrm{inc}}\longrightarrow D_{p,q} \stackrel{[\, \,]}\longrightarrow N_{p,q} \longrightarrow 0.$$
As $E$ is flat, we get another short exact sequence, $$0 \longrightarrow D_{p+1,q-1} \otimes_{A} E \stackrel{\inc \otimes \id}\longrightarrow D_{p,q} \otimes_{A} E \stackrel{[\, \,] \otimes \id}\longrightarrow N_{p,q} \otimes_{A} E \longrightarrow 0$$ but  by definition we also have $$0 \longrightarrow D_{p+1,q-1} \otimes_{A} E \stackrel{\inc \otimes \id}\longrightarrow D_{p,q} \otimes_{A} E \stackrel{[\, \,]_{p,q}}\longrightarrow M_{p,q}  \longrightarrow 0.$$and the result follows from Lemma \ref{raneem}.\quad $\square$

\begin{propos}\label{g2}
If $E$ is a flat left A module, and $\iota : B \longrightarrow A$ is a fibering in the sense of definition \ref{b61}, then $$\Omega^{p} B \otimes_{B} N_{0,q} \otimes_{A} E \cong M_{p,q} $$ via the map
$$ \xi \otimes [x] \otimes e \longmapsto [\iota_{*} \xi \wedge x \otimes e ]_{p,q}.$$
\textbf{Proof:}  Definition  \ref{b61} gives an isomorphism 
$$\Omega^{p} B \otimes_{B} N_{0,q} \longrightarrow N_{p,q}$$
by $ \xi \otimes [x] \longmapsto [\iota_{*} \xi \wedge x].$ 
Now use Proposition \ref{aa66}. \quad  $\square$ 
\end{propos}

We now return to the problem of calculating the cohomology of the sequence (\ref{b60}). Take 
$\xi \otimes [x] \otimes e \in \Omega^{p} B \otimes_{B} N_{0,q} \otimes_{A}E$
(for $x \in \Omega^{q}A)$)
 which maps to $[\iota_{*} \xi \wedge x \otimes e] \in M_{p,q}$, and apply the differential $\nabla^{[p+q]}$ to it to get 
\begin{eqnarray}
&&d(\iota_{*} \xi \wedge x) \otimes e + (-1)^{p+q}\ \iota_{*} \xi \wedge x \wedge \nabla e \cr
&=&  \iota_{*} \extd\xi \wedge x \otimes e + (-1)^{p}\ \iota_{*} \xi \wedge \extd x \otimes e + (-1)^{p+q}\ \iota_{*} \xi \wedge x \wedge \nabla e\ .
\end{eqnarray}
 But $\extd\xi \in \Omega^{p+1} B$, and $$M_{p,q+1} = \frac{\iota_{*} \Omega^{p}B \wedge \Omega^{q+1} A \otimes_{A}E}{\iota_{*} \Omega^{p+1}B \wedge \Omega^{q} A \otimes_{A}E}\ ,$$
 so the first term vanishes on applying $[\quad]_{p,q+1}$. Then 
\begin{eqnarray}\label{g1}
\extd[\iota_{*} \xi \wedge x \otimes e]_{p,q} = (-1)^{p}[\iota_{*} \xi \wedge (\extd x \otimes e +(-1)^{q}x \wedge \nabla e)]_{p,q+1}
\end{eqnarray}
Then, using Proposition \ref{g2}, we have an isomorphism 
\begin{eqnarray}\label{b1}
 \Omega^{p} B \otimes_{B} M_{0,q} \cong  M_{p,q} \ ,\quad \xi \otimes [y]_{0,q} \longmapsto [\iota_{*} \xi \wedge y]_{p,q}\ ,
\end{eqnarray}
 and using this isomorphism,  $\extd$ on $M_{p,q}$ can be written as (see \ref{g1})
 \begin{eqnarray}\label{g4}
\extd( \xi \otimes [y]_{0,q} ) = (-1)^{p} \xi \otimes [\nabla^{[q]} y]_{0,q+1}
\end{eqnarray} 
 where $y \in \Omega^{q} A \otimes_{A} E$. From (\ref{g4}) we see that we should study $[\nabla^{[q]}] : M_{0,q} \longrightarrow M_{0,q+1}$, defined by $[y]_{0,q} \longmapsto [\nabla^{[q]} y]_{0,q+1}$.
 
\begin{propos}  \label{cvauiuy}
$[\nabla^{[q]} ] : M_{0,q} \longrightarrow M_{0,q+1}$ is a left $B$ module map. The module structure is $b . [\eta \otimes e] = [i(b) \eta \otimes e ]$, for $b \in B$ and $\eta \otimes e  \in \Omega^{q} A \otimes_{A} E$. 
\end{propos}
\textbf{Proof:} First, 
\begin{eqnarray*}
[\nabla^{[q]}](b .[\eta \otimes e]_{0,q}) &=& [\extd(i(b)\eta ) \otimes e+(-1)^{q}i(b)\eta \wedge \nabla e]_{0,q+1}\cr
&=& [\iota_{*}(\extd b)\wedge \eta \otimes e+i(b).\extd\eta \otimes e+(-1)^{q}i(b) \eta \wedge \nabla e]_{0,q+1}
\end{eqnarray*}
Now $$\iota_{*}(\extd b)\wedge \eta \otimes e \in \iota_{*} \Omega^{1} B \wedge \Omega^{q}A \otimes_{A}E$$
so $[\iota_{*}(\extd b)\wedge \eta \otimes e]_{0,q+1} = 0$ in $M_{0,q+1}$. Then
\begin{eqnarray*}
[\nabla^{[q]}](b . [\eta \otimes e]_{0,q})& =& [i(b).\extd\eta \otimes e+(-1)^{q}i(b) \eta \wedge \nabla e]_{0,q+1} \cr
&=&b .[\extd\eta \otimes e + (-1)^{q} \eta \wedge\nabla e]_{0,q+1}. \quad \square
\end{eqnarray*}

\begin{propos}\label{bob2}
If  $\Omega^p B$ is flat as a right B module, the cohomology of the cochain complex $$\cdots M_{p,q-1} \stackrel{\extd} \longrightarrow M_{p,q} \stackrel{\extd} \longrightarrow M_{p,q+1} \stackrel{\extd} \longrightarrow  \cdots$$ is given by $\Omega^{p}B \otimes_{B}\hat{H}_{q}$, where $\hat{H}_{q}$ is defined as the cohomology of the cochain complex
$$\cdots \stackrel{[\nabla^{[q-1]}]} \longrightarrow M_{0,q} \stackrel{[\nabla^{[q]}]} \longrightarrow M_{0,q+1} \stackrel{[\nabla^{[q+1]}]} \longrightarrow  \cdots.$$
If we write $\left\langle  \quad \right\rangle_{p,q}$ for the equivalence class in the cohomology of $M_{p,q}$, this isomorphism is given by, for $\xi \in \Omega^{p} B$ and $x \in \Omega^{q} A \otimes_{A} E$, 
\begin{eqnarray}\label{u2}
\left\langle \iota_{*} \xi \wedge x \right\rangle_{p,q} \longrightarrow \xi \otimes \left\langle  x \right\rangle_{0,q}\ .
\end{eqnarray}
\end{propos}
\textbf{Proof:} To calculate the cohomology, we need to find 
${Z}_{p,q} = \im \, \extd : M_{p,q-1} \to M_{p,q}$ and 
${K}_{p,q} = \ker \, \extd : M_{p,q} \to M_{p,q+1}$. 
As we know from Proposition \ref{cvauiuy} that $\extd=[\nabla^{[q]}]: M_{0,q} \longrightarrow M_{0,q+1}$ is a left B module map, we have an exact sequence of left $B$ modules, where the first map is inclusion,
\begin{eqnarray}\label{where}
0 \longrightarrow K_{0,q} \stackrel{\mathrm{inc}} \longrightarrow M_{0,q} \stackrel{\extd} \longrightarrow Z_{0,q+1} \longrightarrow 0\ .
\end{eqnarray}
Since $\Omega^{P}B$ is flat as a right $B$ module, we have another  exact sequence,
\begin{eqnarray}\label{b2}
0 \longrightarrow \Omega^{p}B\otimes_{B} K_{0,q} \stackrel{\id \otimes \inc} \longrightarrow \Omega^{p}B\otimes_{B} M_{0,q} \stackrel{\id \otimes \extd} \longrightarrow \Omega^{p}B\otimes_{B} Z_{0,q+1} \longrightarrow 0\ .
\end{eqnarray} 
Now refer to the isomorphism given in (\ref{b1}), and then by (\ref{g4}) the last map $\id \otimes \extd$ is $(-1)^{p} \extd$ on $M_{p,q}$, so  
 ${Z}_{p,q} = \Omega^{p} B \otimes_{B} Z_{0,q}$ and ${K}_{p,q} = \Omega^{p} B \otimes_{B} K_{0,q}$.

From the definition of $\hat{H}_{q}$ we have another short exact sequence,
\begin{eqnarray*}
0 \longrightarrow Z_{0,q} \stackrel{\inc}\longrightarrow K_{0,q} \longrightarrow \hat{H}_{q} \longrightarrow 0\ ,
\end{eqnarray*} 
and  applying $\Omega^{p}B \otimes_{B}$ gives, as $\Omega^{p}B$ is flat as a right B module, 
\begin{eqnarray}\label{g6}
0 \longrightarrow \Omega^{p} B \otimes_{B} Z_{0,q} \stackrel{\id \otimes \inc}\longrightarrow \Omega^{p} B \otimes_{B} K_{0,q} \longrightarrow \Omega^{p} B \otimes_{B} \hat{H}_{q} \longrightarrow 0\ .
\end{eqnarray}
We deduce that the cohomology of $M_{p,q}$ is isomorphic to $\Omega^{p}B \otimes_{B} \hat{H}_{q}$. \quad
 $\square$

\subsection{The second page of the spectral sequence}
Now we move to the second page  of the spectral sequence, in which we take the cohomology of the previous cohomology, i.e. the cohomology of  
$$\extd : \mathrm{cohomology} \,  (M_{p,q} ) \longrightarrow \mathrm{cohomology} \, (M_{p+1,q}).$$
By the isomorphism discussed in Proposition \ref{bob2}, we can view this as 
\begin{eqnarray}
\extd : \Omega^{p}B \otimes_{B} \hat{H}_{q}\longrightarrow \Omega^{p+1}B \otimes_{B} \hat{H}_{q}\ .
\end{eqnarray}

\begin{propos}\label{b6}
The differential $\extd$ gives a left covariant derivative $$\nabla_{q} : \hat{H}_{q} \longrightarrow \Omega^{1}B \otimes_{B} \hat{H}_{q}.$$ If $\left\langle \xi \otimes e\right\rangle_{0,q} \in \hat{H}_{q}$, this is given by using the isomorphism (\ref{u2}) as
$$\left\langle \xi \otimes e\right\rangle_{0,q} \longmapsto \eta \otimes \left\langle \omega \otimes f \right\rangle_{0,q}\ ,$$
 where 
 $$\extd\xi \otimes e + (-1)^{q}\xi \wedge \nabla e = \iota_{*} \eta \wedge \omega \otimes f
 \in \iota_*\Omega^1 B\wedge\Omega^q A\tens_A E\ .
 $$  
 \end{propos}
\textbf{Proof}: Take $\left\langle x\right\rangle_{0,q} \in \hat{H}_{q}$, where $x \in K_{0,q} = \ker \, \extd : M_{0,q} \to M_{0,q+1}$, and suppose $x = \xi \otimes e$, where $\xi \in \Omega^{q}A$ and $e \in E$ (summation implicit). As $x \in K_{0,q}$ we have 
$$[\extd x]_{0,q+1} = [\extd\xi \otimes e + (-1)^{q}\xi \wedge \nabla e]_{0,q+1} = 0$$
 in $M_{0,q+1}$, so 
 $$\extd\xi \otimes e + (-1)^{q}\xi \wedge \nabla e \in \iota_{*} \Omega^{1}B \wedge \Omega^{q} A \otimes_{A} E.$$ We write (summation implicit), for $\eta \in \Omega^{1} B$,
 $\omega \in \Omega^{1} A$ and $f \in E$,
\begin{eqnarray}\label{b3}
\extd\xi \otimes e +(-1)^{q} \xi \wedge \nabla e = \iota_{*}\eta \wedge \omega \otimes f\ .
\end{eqnarray}
Under the isomorphism (\ref{b1}), this corresponds to $\eta \otimes [\omega \otimes f]_{q} \in \Omega^{1}B \otimes_{B} M_{0,q}$.
As the curvature of $E$ vanishes, we have from applying $\nabla^{[q+1]}$ to (\ref{b3}), 
\begin{eqnarray}\label{b4}
\iota_{*} \extd\eta \wedge \omega \otimes f - \iota_{*} \eta \wedge \extd\omega \otimes f + (-1)^{q+1} \iota_{*} \eta \wedge \omega \wedge \nabla f = 0 .
\end{eqnarray}
We take this as an element of $M_{1,q+1}$, so we apply $[\quad ] _{1,q+1}$ to (\ref{b4}). Then as the denominator of $M_{1,q+1}$ is $\iota_{*} \Omega^{2}B \wedge \Omega^{q} A \otimes_{A} E$, we see that the first term of (\ref{b4}) vanishes on taking the quotient, giving
$$- [ \iota_{*} \eta \wedge (\extd\omega \otimes f +(-1)^{q} \omega \wedge \nabla f )] _{1,q+1} = 0.$$ Under the isomorphism (\ref{b1}) this corresponds to 
\begin{eqnarray}\label{b5}
- \eta \otimes_{B} [\extd\omega \otimes f + (-1)^{q} \omega \wedge \nabla f ]_{0,q+1} = 0.
\end{eqnarray}
This means that $$\eta \otimes [ \omega \otimes f ]_{0,q} \in \Omega^{1}B \otimes_{B} M_{0,q}$$
is in the kernel of the map $\id \otimes \extd$ in (\ref{b2}), and as  (\ref{b2}) is an  exact sequence  we have $$\eta \otimes [\omega \otimes f]_{0,q} \in \Omega^{1} B \otimes_{B} K_{0,q},$$ so we can see take the cohomology class to get $$\eta \otimes \left\langle  \omega \otimes f \right\rangle_{0,q} \in \Omega^{1}B \otimes_{B} \hat{H}_{q}.$$
This completes showing that $\nabla_{q}$ exists, but we also need to show that it is a left covariant derivative. For $b \in B$, we calculate $\nabla_{q}(b. \xi \otimes e )$ to get $$\extd(b.\xi) \otimes e + (-1)^{q} b. \xi \wedge \nabla e = \extd b \wedge \xi \otimes e + b. (\extd\xi \otimes e + (-1)^{q} \xi \wedge \nabla e),$$
so we get $$\nabla_{q}\left\langle b. \xi \otimes e\right\rangle _{0,q} = \extd b \otimes \left\langle \xi \otimes e\right\rangle_{0,q} + b. \nabla_{q}\left\langle \xi \otimes e \right\rangle _{0,q}\ .    \quad \square$$

\begin{propos}\label{b9}
The curvature $R_q$ of the covariant derivative $\nabla_{q}$ in Proposition \ref{b6} is zero.
\end{propos}
\textbf{Proof:} Using the notation of Proposition \ref{b6}, equation (\ref{b3}) $$\nabla_{q}\left\langle \xi \otimes e\right\rangle_{0,q} = \eta \otimes \left\langle \omega \otimes f\right\rangle_{0,q}.$$ 
If we apply $\nabla_{q}^{[1]}$ (see Definition \ref{key49}) to this, we get 
\begin{eqnarray}\label{u4}
R_{q}\left\langle \xi \otimes e\right\rangle _{0,q} = \extd\eta \otimes \left\langle \omega \otimes f\right\rangle_{0,q} - \eta \wedge \nabla_{q} \left\langle \omega \otimes f\right\rangle _{0,q}.
\end{eqnarray}
To find $\nabla_{q}\left\langle w \otimes f\right\rangle _{0,q}$, referring to the proof of Proposition \ref{b6}, formula (\ref{b5}), we have 
$$\eta \otimes_{B}(\extd\omega \otimes f + (-1)^{q} \omega \wedge \nabla f ) \in \Omega^{1}B \otimes_{B} (\iota_{*} \Omega^{1}B \wedge \Omega^{q}A \otimes_{A}E).$$
This comes for tensoring the exact sequence
$$0 \longrightarrow \iota_{*} \Omega^{1}B \wedge \Omega^{q}A \otimes_{A}E \longrightarrow \Omega^{q+1}A \otimes_{A} E \stackrel{[\, \,]_{0,q+1}}\longrightarrow M_{0,q+1}\longrightarrow 0$$ on the left by $\Omega^{1}B$, and using that $\Omega^{1}B$ is a flat right module. Now write (summation implicit),
\begin{eqnarray}\label{u5}
\eta \otimes (\extd\omega \otimes f +(-1)^{q} \omega \wedge \nabla f) = \eta^{\prime} \otimes (\iota_{*} \kappa \wedge \zeta \otimes g) 
\end{eqnarray} 
for $\eta^{\prime}$,$\kappa \in \Omega^{1} B$,  $\zeta \in \Omega^{q}A$ and $g \in E$.
Then, from Proposition \ref{b6},
$$\eta \wedge \nabla_{q}\left\langle \omega \otimes f\right\rangle _{0,q} = \eta^{\prime} \wedge \kappa \otimes \left\langle \zeta \otimes g\right\rangle _{0,q}$$ so from (\ref{u4}),
\begin{eqnarray}\label{arwa}
R_{q}\left\langle \xi \otimes e\right\rangle _{0,q} = \extd\eta \otimes \left\langle \omega \otimes f\right\rangle_{0,q} - \eta^{\prime} \wedge \kappa \otimes \left\langle  \zeta \otimes g\right\rangle_{0,q}\ .
\end{eqnarray} 
Now (\ref{u5}) implies that $$\iota_{*}\eta \wedge (\extd \omega \otimes f + (-1)^{q} \omega \wedge \nabla f) = \iota_{*} \eta^{\prime} \wedge \iota_{* } \kappa \wedge \zeta \otimes g\ ,$$ 
and substituting this into (\ref{b4}) gives 
$$\iota_{*} \extd\eta \wedge \omega \otimes f - \iota_{*} \eta^{\prime} \wedge \iota_{*} \kappa \wedge \zeta \otimes g = 0\ ,$$
so on taking equivalence classes in $M_{2,q}$ we find, using the isomorphism (\ref{b1}),
$$\extd\eta \otimes [\omega \otimes f]_{0,q} - \eta^{\prime} \wedge \kappa \otimes [\zeta \otimes g]_{0,q} = 0\ ,$$ 
and this shows that $R_q=0$ by (\ref{arwa}).   $\square$

\begin{theorem}
Given:\\
1) a map $\iota : B \longrightarrow A$ which is a differential fibration (see definition \ref{b61}),\\
2) a flat left A module $E$, with a zero-curvature left covariant derivative $\nabla_{E} : E \to \Omega^{1}A \otimes_{A} E$,\\
3) each $\Omega^{p}B$ is flat as a right B module,\\
then there is a spectral sequence converging to $ H^{*}(A, E, \nabla_{E})$ with second page $ H^{*}(B, \hat{H}_{q}, \nabla_{q})$ where $\hat{H}_{q}$ is defined as the cohomology of the cochain complex$$\cdots \stackrel{\extd} \longrightarrow M_{0,q} \stackrel{\extd} \longrightarrow M_{0,q+1} \stackrel{\extd} \longrightarrow  \cdots$$
where 
\begin{eqnarray*}
M_{0,q} &=& \frac{\Omega^{q}A \otimes_{A} E}{\iota_{*} \Omega^{1} B \wedge \Omega^{q-1} A \otimes_{A}E}\ ,\cr
\extd[x \otimes e]_{0,q} &=& [\extd x \otimes e + (-1)^{q} x \wedge \nabla_{E} e]_{0,q+1}\ .
\end{eqnarray*}
The zero curvature left covariant derivative $\nabla_{q}: \hat{H}_{q} \to \Omega^{1}B \otimes_{B} \hat{H}_{q}$ is as defined in Proposition \ref{b6}.
\end{theorem} 
\textbf{Proof:} The first part of the proof is given in Proposition \ref{bob2}. Now we need to calculate the cohomology of $$\extd : \Omega^{p}B \otimes_{B} \hat{H}_{q} \longrightarrow \Omega^{p+1}B \otimes_{B}\hat{H}_{q}$$ This is given for $\xi \otimes \left\langle  \eta \otimes e\right\rangle _{0,q}$ (for $\xi \in \Omega^{p}B$, $\eta \in \Omega^{q} A$ and $e \in E$) as follows: this element corresponds to $\iota_{*}\xi \wedge \eta \otimes e$, and applying $\extd$ to this gives $$\iota_{*} \extd\xi \wedge \eta \otimes e + (-1)^{p} \iota_{*} \xi \wedge \extd\eta \otimes e +(-1)^{p+q} \iota_{*} \xi \wedge \eta \wedge \nabla e.$$
But we have calculated the effect of $\extd$ on $\hat{H}_{q}$ in Proposition \ref{b6}, so we get 
$$\extd(\xi \otimes \left\langle \eta \otimes e \right\rangle _{0,q} ) = \extd\xi \otimes  \left\langle \eta \otimes e \right\rangle _{0,q} +(-1)^{p}\xi \wedge \nabla_{q}  \left\langle \eta \otimes e \right\rangle _{0,q}.$$ 
The covariant derivative $\nabla_{q} $ has zero curvature by Proposition \ref{b9}.  \quad $\square$

\section{Example: A fibration with fiber the noncommutative torus}\label{se1}
As discussed in the Introduction, the idea for this example came from \cite{25,26}.

\subsection{The Heisenberg group}  \label{vcadgsh}
The  Heisenberg group $H$ is defined to be following subgroup of $M_{3}\mathbb{(\mathbb{Z})}$ under multiplication. 
 \[ \Big\{  \left(
\begin{array}{ccc}
1&n&k\\
0&1&m\\
0&0&1
\end{array}
\right)
  : n , m , k \in \mathbb{Z} \Big\} \]
  We can take generators $u, v , w$ for the group,   where $w$ is central and there is one more relation $uv=wvu$.  These generators correspond to the matrices
  \[ u=  \left(
\begin{array}{ccc}
1&1&0\\
0&1&0\\
0&0&1
\end{array}
\right)  ,\quad
 v=  \left(
\begin{array}{ccc}
1&0&0\\
0&1&1\\
0&0&1
\end{array}\right)\ ,\quad
 w=  \left(
\begin{array}{ccc}
1&0&1\\
0&1&0\\
0&0&1
\end{array}\right)\
\]
There is an isomorphism $\theta : H \longrightarrow H $, for every matrix
 \[\left(
\begin{array}{cc}
a&c\\
b&d
\end{array}\right)\
  \in SL(2,\mathbb{Z}) ,
\]
given by $\theta(u)=u^{a}\,v^{b}$,    $\theta(v)=u^{c}\,v^{d}$,      $\theta(w)=w$.
The group algebra $\mathbb{C}H$  of $H$ can be made into a star algebra by setting $x^*=x^{-1}$ for all $x \in \{ u , v , w \}$.

\subsection{A differential calculus on the Heisenberg group}
There is a differential calculus on the group algebra $\mathbb{C}H$  of $H$. 
It is bicovariant, as set down by Woronowicz in \cite{worondiff}. 

For a generator $ x \in \{u , v , w\}$,  we write $e^x = x^{-1}.\extd x$, a left invariant element of $\Omega^{1}\mathbb{C}H$. 
We suppose that $\Omega^{1}\mathbb{C}H$ is free as left $\mathbb{C}H$ module, with generators $\{ e^{u} , e^{v} , e^{w} \}$. This means that every element of      $\Omega^{1}\mathbb{C}H$ can be written uniquely as $a^{u} . e^{u} + a^{v} . e^{v} + a^{w} . e^{w}$, for $a^{u} , a^{v} , a^{w} \in \mathbb{C}H$. We have the following relations on $\Omega^{1}\mathbb{C}H$, for all $x \in \{ u , v , w \}$:\\
$x .e^x  = e^x .x $   \\
$x .e^w  = e^w . x $\\
$w .e^x  = e^x .w $\\
$ u^{-n} .e^v. u^n =  e^v - \frac{n}{2}\,e^w $ \\
$ v^{-n}. e^u. v^n =  e^v + \frac{n}{2}\,e^w $\\
Further the map $\theta$ in subsection \ref{vcadgsh}
extends to a map of 1-forms given by 
$     \theta(e^{w}) = e^{w}    $   \\
$      \theta(e^{u}) = a .e^{u} + b. e^{v} + \frac{ab}{2}. e^{w}   $   \\
$     \theta(e^{v}) = c .e^{u} + d. e^{v} + \frac{cd}{2}. e^{w}    $   \\
Checking the braiding given by Woronowicz shows that, for $x,y \in \{u , v , w\}$,\\
$\extd e^x=0$ \\
$e^x\wedge e^y=-e^y\wedge e^x$\ .\\
The star operation extends to the differential calculus, with $(e^x)^*=-e^x$.

\subsection{The differential fibration}\label{b10}
If we take $z$ to be the identity function :$S^{1} \rightarrow \mathbb{C}$, the map sending $z^n$ to $w^n$ gives an algebra map $\iota:C(S^1)\to \mathbb{C}H$. It is also a star algebra map, with the usual star structure $z^*=z^{-1}$ on $C(S^1)$. 

The differential structure of the `fiber algebra' $F$ is 
\begin{eqnarray}\label{fib1}
\Omega^n F = \frac{\Omega^n \mathbb{C}H}{\iota_*\Omega^1 C(S^1) \wedge \Omega^{n-1}\mathbb{C}H}\ ,
\end{eqnarray}
i.e.\ we put $ \extd w=0$  in $\Omega^n F$ (i.e. put $e^{w} = 0$). This is because in (\ref{fib1}) we divide by everything of the form $e^{w} \wedge \xi $. To see that this gives a fibration, we note that a linear basis for the left invariant n-forms is as follows:\\
$\Omega^{1}A$:  \quad $e^{u},  e^{v},  e^{w}$\\
$\Omega^{2}A$: \quad  $e^{u} \wedge e^{v}$, $e^{w} \wedge e^{u}$ , $e^{w} \wedge e^{v}$\\
$\Omega^{3}A$: \quad  $e^{v} \wedge e^{u} \wedge e^{w}$\\
Then the $N_{n,m}$ (see (\ref{cvhgsuv})) are, where $\<...\>$ denotes the module generated by, and all others are zero:\\
$N_{0,0}  = 1$, \quad  $N_{1,0} = \left\langle e^{w} \right\rangle $, \quad $N_{m,0} = 0$ ,  $m > 1$ \\
$N_{0,1}  = \frac{\left\langle  e^{u}, e^{v}, e^{w}\right\rangle }{\left\langle  e^{w} \right\rangle } = \left\langle e^{u}, e^{v} \right\rangle $\\
$N_{0,2}  = \frac{\left\langle  e^{u} \wedge e^{v}, e^{w} \wedge e^{u}, e^{w} \wedge e^{v} \right\rangle }{\left\langle e^{w} \wedge e^{u}, e^{w} \wedge e^{v} \right\rangle } = \left\langle e^{u} \wedge e^{v} \right\rangle $\\
$N_{0,3}  = \frac{\left\langle  e^{w} \wedge e^{u} \wedge e^{v}\right\rangle }{\left\langle   e^{w} \wedge e^{u} \wedge e^{v} \right\rangle } = 0$\\
$N_{0,n} = 0 $ \quad $n \geq 4$\\
$N_{1,1} =  \frac{e^{w} \wedge \left\langle e^{w}, e^{u}, e^{v} \right\rangle }{\left\langle 0 \right\rangle } = \left\langle e^{w} \wedge e^{u} , e^{w} \wedge e^{v} \right\rangle $\\
$N_{1,2} = \frac{e^{w} \wedge \left\langle e^{w} \wedge e^{u}, e^{w} \wedge e^{v}, e^{u} \wedge e^{v} \right\rangle }{\left\langle 0 \right\rangle } = \left\langle e^{w} \wedge e^{u} \wedge e^{v} \right\rangle$ \\
Then the following map  is one-to-one and onto, 
$$\Omega^{1}C(S^{1}) \otimes_{C(S^{1})} N_{0,n}   \longrightarrow N_{1,n}$$
 giving a differential fibration in the sense of Definition \ref{b61}.

As was done in \cite{25}, we note that this map does have a fiber in quite a classical sense.
The algebra $C(S^1)$ is commutative, and if we take $q\in S^1$, the fiber algebra
corresponding to $q$ is given by substituting $w\mapsto q$ in the algebra relations. We get unitary generators $u ,v$ and a relation $u\, v = q\, v\, u$ for a complex number $q$ of norm $1$.  But this is exactly the noncommutative torus $\mathbb{T}^{2}_{q}$. The map $\theta$ on the total algebra
$\mathbb{C}H$ is the identity on the base algebra $C(S^1)$, so it acts on each fiber.


\begin{thebibliography}{ggghhh}


 \bibitem{48}
	  Andrea F., Dabrowski L.\ \& Landi G.,
	  \emph{ The Isospectral Dirac Operator on the 4-dimensional Orthogonal Quantum Sphere.}
	 Commun. Math. Phys. 279, 77-116 (2008).



	 \bibitem{three}
Beggs E.J.\ \& Brzezi\'nski T.,
	\emph{The Serre spectral sequence of noncommutative fibration for de Rham cohomology}. Acta Mathematica 195 (2005), p155-196 


      \bibitem{44}
	  \ Berrick A.J.\ \&  Keating M.E.,
	  \emph{ An Introduction to Ring and Modules with K-theory in view.}
	 C.U.P. 2000.
	 
\bibitem{BWquant93} Biedenharn L.C.\ \&  Lohe M.A.,
	 \emph{An extension of the Borel-Weil construction to the quantum group $U_q(n)$},
	 Comm.\ Math.\ Phys., Volume 146, Number 3, 483-504 (1992)
	 
\bibitem{BraMa}
{Brain S.J.\ \& Majid S.,}
\emph{Quantisation of twistor theory by cocycle twist}, 
Commun.\ Math.\ Phys., vol. 284, no.\ 3, pp.\ 713-774, 2008

	\bibitem{28}
	   Bredon G.E.,
	  \emph{Sheaf Theory,}
	1967 by McGraw-Hill.
	
	 \bibitem{17}
	  Bresser K., Muller-Hoissen F., Dimakis A.\ \&  Sitarz A.,
	  \emph{Noncommutative geometry of finite groups.}
	  J. of Physics A (Math. and General), 29 :2705- 2735, 1996. 
 
\bibitem{ConnesIHES}
{Connes A.,} \emph{Non-commutative differential geometry},
Publications Math\'ematiques de L'IH\'ES,
 Vol.\ 62, No.\ 1, Pages 41-144, 1985

 \bibitem{25}
Echterhoff S., Nest R.\ \& Oyono-Oyono H.,
	  \emph{Principal noncommutative torus bundles,}
Proc.\ London Math.\ Soc.\ (2009) 99 (1): 1-31.

\bibitem{26}
Echterhoff S., Nest R.\ \& Oyono-Oyono H.,
	  \emph{Fibration with noncommutative fibers,}
	J.\ Noncommut.\ Geom.\ 3 (2009), no.\ 3, 377-417

\bibitem{Ma:rief}
{Majid S.,} \emph{Riemannian geometry of quantum groups and
finite groups with nonuniversal differentials},   Commun.\ Math.\ Phys.\
225 (2002), 131-170.

 \bibitem{MasThesis}
{Masmali I.,} \emph{Hopf algebra and noncommutative differential structures}, Ph.D.\ thesis, Swansea University, 2010

	   \bibitem{11}
	  McCleary J.,
	  \emph{A User's Guide to Spectral Sequences.}
	  2nd ed., Cambridge University Press, Cambridge (2001).

\bibitem{MulQuant} Mulvey C.J. \& Pelletier J.W., \emph{On the Quantisation of Spaces},  J.\ Pure Appl.\ Algebra 175, no.\ 1-3, 289-325  (2002)



	   \bibitem{quillModel}
	  Quillen D.G., \emph{Homotopical algebra}, Lecture Notes in Mathematics, Springer-Verlag, no.\ 43 (1967).
 
\bibitem{worondiff}
{Woronowicz S.L.,}
\emph{Differential calculus on compact matrix pseudogroups (quantum groups)},
 Comm.\ Math.\ Phys.\  122  (1989),  no.\ 1, 125--170. 



\end{thebibliography}
\end{document}